\documentclass[12pt]{article}
\usepackage{amsmath, amsfonts, amsthm, amssymb, epsfig, textcomp}

\newtheorem{theorem}{Theorem}
\newtheorem{lemma}{Lemma}

\begin{document}
\title{{\bf The Rational Number} $\mathbf{\frac{n}{p}}$ {\bf as a sum of two
  unit fractions}}
\author{Konstantine Zelator\\
Department of Mathematics, Statistics, and Computer Science\\
212 Ben Franklin Hall\\
Bloomsburg University of Pennsylvania\\
400 East 2nd Street\\
Bloomsburg, PA  17815\\
USA\\
and\\
P.O. Box 4280\\
Pittsburgh, PA  15203\\
e-mails: kzelator@bloomu.edu\\ 
\hspace{1.5in} konstantine\underline{\ }zelator@yahoo.com}

\maketitle

\section{Introduction}  In a 2011 paper in the journal {\it Asian Journal of
  Algebra} (see \cite{1}), the authors consider, among other equations, the
  diophantine equations

$$2xy= n(x+y)\ \ {\rm and}\ \ 3xy=n(x+y).
$$

For the first equation, with $n$ being an odd positive integer, they
give the solution (in positive integers $x$ and $y$) $\dfrac{n+1}{2} = x =
\dfrac{(n-1)}{2} + 1$, $y=n\left(\dfrac{(n-1)}{2} + 1 \right)  = n \left(
\dfrac{n+1}{2}\right)$.  

For the second equation, with $n \equiv 2({\rm mod}3)$, they present the
particular solution,

$$
\dfrac{n+1}{3} = x = \dfrac{(n-2)}{3} + 1,\ \ y = n\left(\dfrac{(n-2)}{3} + 1
\right) = n \left(\dfrac{n+1}{3}\right).
$$

If in the above equations we assume $n$ to be prime, then these two equations
become special cases of the diophantine equation,  $nxy = p(x+y)$, with $p$
being a prime and $n$ a positive integer with $n \geq 2$.

This two-variable symmetric diophantine equation is the subject matter of this
article; with the added condition that the integer $n$ is not divisible by the
prime $p$.  Observe that this equation can be written equivalently in fraction
form:

$$
\dfrac{n}{p} = \dfrac{1}{x} + \dfrac{1}{y}.
$$

\noindent This problem then can be approached from the point of view of
decomposing a positive rational number into a sum of two unit fractions (i.e.,
two rational numbers whose numerators are equal to $1$).  The ancient
Egyptians left behind an entire body of work involving the decomposition of a
given fraction into a sum of two or more unit fractions.  They did so by
creating tables containing the decomposition of specific fractions into sums
of unit fractions.  An excellent source on the subject of the work of the
ancient Egyptians on unit fractions is the book by David M. Burton, ``The
History of Mathematics, An Introduction'' (see \cite{2}).  Note that thanks to
the identity $\dfrac{1}{k} = \dfrac{1}{k+1} + \dfrac{1}{k(k+1)}$, a unit
fraction can always be written as a sum of two unit fractions.

We state our theorem.

\begin{theorem}  Let $p$ be a prime, $n$ a positive integer, $n \geq 2$.
  Also, assume that gcd$(p,n)=1$ (equivalently, $n$ is not divisible by
  $p$).  Consider the two-variable symmetric diophantine equation, 

\begin{equation} nxy = p (x+y) \label{E1}
\end{equation}

\noindent with the two variables $x$ and $y$ taking values from the set
${\mathbb Z}^+$ of positive integers.  Then, 

\begin{enumerate}
\item[(i)] If $n =2$ and $p \geq 3$, equation (\ref{E1}) has exactly three
  distinct solutions, the following positive integer pairs:

$$
(x,y) = (p,p),\ (x,y) = \left(p\left(\dfrac{p+1}{2}\right),\ \dfrac{p+1}{2}
  \right),
$$

\noindent and its symmetric counterpart 

$$
(x,y) = \left( \dfrac{p+1}{2} , \ p\left(\dfrac{p+1}{2}\right)\right).
$$

\item[(ii)]  If $n \geq 3$, and $n$ is a divisor of $p+1$.  Then equation
  (\ref{E1}) has exactly two distinct solutions:

$$
(x,y) = \left(p\left(\dfrac{p+1}{n} \right),\ \dfrac{p+1}{n} \right) \ {\rm
  and} \ (x,y) = \left( \dfrac{p+1}{n},\ p \left(\dfrac{p
  +1}{n}\right)\right).$$

\item[(iii)]  If $n$ is not a divisor of $p+1$, Equation (\ref{E1}) has no
  solution.
\end{enumerate}
\end{theorem}

\section{A lemma from number theory}

The following lemma, commonly referred to as Euclid's lemma, is of great
significance in number theory.

\begin{lemma}  (Euclid's lemma):  Suppose that $a,b,c$ are positive integers
  such that $a$ is a divisor of the product $bc$; and gcd$(a,b)=1$ (i.e., $a$
  and $b$ are relatively prime),  then $a$ must be a divisor of $c$.
\end{lemma}

Typically, this lemma and its proof can be found in an introductory number
theory book.  For example, see reference \cite{3}.

\section{ Proof of Theorem 1}

First we show that the positive integer pairs listed in Theorem 1 are indeed
solutions to Equation (\ref{E1}).

If $n = 2$ and $p\geq 3$, then for $(x,y)=(p,p)$, a straightforward
calculation shows both sides of (\ref{E1}) are equal to $2p^2$;
and for $(x,y) = \left( \dfrac{p(p+1)}{2}, \dfrac{p+1}{2} \right)$, a
calculation shows that both sides of (\ref{E1}) are equal to
$\dfrac{p(p+1)^2}{2}$.  

If $n \geq 3$ and $n$ is a divisor of $p+1$, then for $(x,y) = \left( p
\left(\dfrac{p+1}{n} \right), \dfrac{p+1}{n} \right)$, a calculation shows
that both sides of equation (\ref{E1}) are equal to $\dfrac{p(p+1)^2}{n}$.

In the second part of this proof, we show that there are no other solutions to
equation (\ref{E1}).  To do so, we will demonstrate that if $(t_1,t_2)$ is a
solution to (\ref{E1}), then it must be one of the solutions listed in Theorem
1.  So, let $(t_1,t_2)$ be a positive integer solution to equation (\ref{E1}).

We have,

\begin{equation}
\left\{ \begin{array}{c} p(t_1 + t_2) = nt_1 t_2\\
\\
t_1, t_2 \in {\mathbb Z}^+ \end{array}\right\} \label{E2}
\end{equation}

\noindent Let $d$ be the greatest common divisor of $t_1$ and $t_2$.  Then 

\begin{equation}
\left\{ \begin{array}{l} t_1 = du_1,\ t_2 = du_2;\\
{\rm for\ relatively \ prime\ positive\ integers}\ u_1\ {\rm and}\ u_2;\\
{\rm gcd}(u_1,u_2) = 1 \end{array} \right\} \label{E3}
\end{equation}

\noindent From (\ref{E2}) and (\ref{E3}) we obtain,

\begin{equation}
p(u_1+u_2) = nd\,u_1u_2 \label{E4}
\end{equation}

\noindent Since the prime $p$ is relatively prime to $n$.  By (\ref{E4}) and
Lemma 1, it follows that $p$ must divide the product $du_1u_2$. Since $p$
is a prime number, it must divide at least one of $d$ and $u_1u_2$.  We
distinguish between two cases:  The case wherein $p$ divides the product
$u_1u_2$; and the case in which $p$ is a divisor of $d$.

\vspace{.15in}

\noindent {\bf Case 1:}  {\it $p$ is a divisor of $u_1u_2$.}

Since $p$ is a prime, and the integers $u_1$ and $u_2$ are relatively prime by
(\ref{E3}), and also in view of the fact that $p$ divides the product $u_1u_2$,
it follows that $p$ must divide exactly one of $u_1,u_2$.  It must divide one
but not the other.  Thus, there are two subcases in Case 1.  

\vspace{.15in}

\noindent Subcase 1a being the one with   $p|u_1$ (i.e., $p$ divides $u_1$);

\noindent Subcase 1b:  $p$ divides $u_2$.

\vspace{.15in}

\noindent But these two subcases are symmetric since equation (\ref{E4}) is
symmetric in $u_1$ and $u_2$.  Thus, without loss of generality, we need only
consider the subcase $p|u_1$.  So we set

\begin{equation} \left( u_1=pv_1,\ v_1\ {\rm a \ positive\ integer}\right)
  \label{E5}
\end{equation}

\noindent Combining (\ref{E5}) with (\ref{E4}) we get,

\begin{equation}
\left\{ \begin{array}{c} pv_1 + u_2 = nd\, v_1u_2\\
\\
{\rm or \ equivalently}, \ u_2 = v_1 \cdot (ndu_2-p) \end{array}\right\}
\label{E6}
\end{equation}

\noindent According to (\ref{E6}), the positive integer $v_1$ is a divisor of
$u_2$.  But, by (\ref{E5}) $v_1$ is also a divisor of $u_1$.  Since $u_1$ and
$u_2$ are relatively prime by (\ref{E3}), it follows that

\begin{equation}
v_1=1 \label{E7}
\end{equation}

\noindent Hence, by (\ref{E7}) and (\ref{E6}), we further obtain,

\begin{equation}
p = u_2 (nd-1) \label{E8}
\end{equation}

\noindent According to (\ref{E8}), $u_2$ is a divisor of $p$, and since $p$ is
a prime it follows that either $u_2=1$ or $u_2 = p$.  If $u_2 =1$, then
(\ref{E8}) yields $p+1 = nd$ which implies that $n$ is a divisor of $p+1$.
Using $d=\dfrac{p+1}{n}, \ v_1=1, u_2 = 1$, we also get $u_1 = p$ (by
(\ref{E5})). So, by (\ref{E3}) we obtain the solution $t_1 =
p\left(\dfrac{p+1}{n} \right), t_2 = \dfrac{p+1}{n}$ (already a verified
solution in the first part of the proof).  Now, if $u_2 = p$ in (\ref{E8}),
then $2 = nd$ which implies either $n = 2$ and $d=1$, or $n=1$ and $d=2$.  But
$n \geq 2$, so the latter possibility is ruled out.  Thus, $u_2 = p, \ n=2$,
and $d=1$.  Also, by (\ref{E7}) we have $v_1=1$ and so $u_1=p$ by (\ref{E5}).

Hence, (\ref{E3}) yields $t_1 = p = t_2$;  $(p,p)$ with $n=2$ being a
solution verified in the first part of the proof.

\vspace{.15in}

\noindent {\bf Case 2:}  {\it $p$ is a divisor of $d$}

We set

\begin{equation} (d=p\delta,\ \delta\ {\rm is \ a\ positive\
    integer})\label{E9}
\end{equation}

\noindent by (\ref{E9}) and (\ref{E4}) we have,

\begin{equation}
u_1 + u_2 = n\delta u_1u_2 \label{E10}
\end{equation}

Clearly, by inspection, we see that equation (\ref{E10}) implies that the
positive integers $u_1$ and $u_2$ must divide each other.  Since they are
relatively prime, it follows that 

\begin{equation}
u_1=u_2 = 1 \label{E11}
\end{equation}

Equations (\ref{E10}) and (\ref{E11}) yield

\begin{equation}
2 = n \delta \label{E12}
\end{equation}

Due to the fact  $n \geq 2$, (\ref{E12}) implies that $n=2$ and $\delta =1$.
So, by (\ref{E11}), (\ref{E9}), and (\ref{E3}), it is clear that (since $d=p$)
$u_1=u_2=p$.  This produces $(u_1,u_2)=(p,p)$, with $n=2$.  An already
verified solution.  The proof is complete. \hfil $\Box$


\begin{thebibliography}{99}
\bibitem[1]{1} Kishan, Hari, Rani, Megha and Agarwal, Smiti, The Diophantine
  Equations of Second and Higher Degree of the Form $3xy=n(x+y)$ and
  $3xyz=n(xy+yz + zx)$, etc., {\it Asian Journal of Algebra} 4(1), (2011), pp. 31-37.

\bibitem[2]{2} Burton, David M., ``The History of Mathematics, An
  Introduction'', Sixth Edition, McGraw Hill, (2007),  p. 40.

\bibitem[3]{3}  Rose, Kenneth H., ``Elementary Number Theory and Its
  Applications'', 5th Edition, Pearson, Addison Wesley, (2005), p. 109.
\end{thebibliography}
\end{document}